\documentclass[12pt]{article}
\usepackage{latexsym,amsfonts,amsthm,amsmath,amscd,amssymb}
\usepackage[dvips]{graphicx}

\newtheorem{lemma}{Lemma}[section]
\newtheorem{thm}[lemma]{Theorem}
\newtheorem{rem}[lemma]{Remark}
\newtheorem{prop}[lemma]{Proposition}
\newtheorem{cor}[lemma]{Corollary}

\newtheorem{defn}[lemma]{Definition}

\reversemarginpar

\newcommand{\dimo}[1]{\vspace{2pt}\noindent\textit{Proof of~\ref{#1}}.\ }
\newcommand{\finedimo}{{\hfill\hbox{$\square$}\vspace{2pt}}}

\newcommand\matT{{\mathbb{T}}}
\newcommand\matS{{\mathbb{S}}}

\newcommand\matC{{\mathbb{C}}}

\newcommand\Sigmatil{{\widetilde\Sigma}}
\newcommand\ntil{{\widetilde n}}
\newcommand\ftil{{\widetilde f}}
\renewcommand{\hbar}{{\overline{h}}}
\newcommand{\ee}{{\rm e}}
\newfont{\Got}{eufm10 scaled 1200}
\newcommand{\permu}{{\hbox{\Got S}}}
\newcommand{\compo}{\,{\scriptstyle\circ}\,}

\newcommand{\mycap} [1] {\caption{\footnotesize{#1}}}

\newcommand\calR{{\mathcal R}}

\begin{document}

\title{On the existence of\\ branched coverings between
surfaces\\ with prescribed branch data, II}

\author{Ekaterina~{\textsc Pervova}\thanks{Partially supported by the INTAS YS
fellowship 03-55-1423} \and\addtocounter{footnote}{5}
Carlo~{\textsc Petronio}\footnote{Supported by the INTAS project
``CalcoMet-GT'' 03-51-3663}}

\maketitle

\begin{abstract}
\noindent If $\Sigmatil\to\Sigma$ is a branched covering between
closed surfaces, there are several easy relations one can
establish between $\chi(\Sigmatil)$, $\chi(\Sigma)$, orientability
of $\Sigma$ and $\Sigmatil$, the total degree, and the local
degrees at the branching points, including the classical
Riemann-Hurwitz formula. These necessary relations have been shown
to be also sufficient for the existence of the covering except
when $\Sigma$ is the sphere $\matS$ (and when $\Sigma$ is the
projective plane, but this case reduces to the case
$\Sigma=\matS$). For $\Sigma=\matS$ many exceptions are known to
occur and the problem is widely open.

Generalizing methods of Bar\'anski, we prove in this paper that
the necessary relations are actually sufficient in a specific but
rather interesting situation. Namely under the assumption that
$\Sigma=\matS$, that there are three branching points, that one of
these branching points has only two preimages with one being a
double point, and either that $\Sigmatil=\matS$ and that the
degree is odd, or that $\Sigmatil$ has genus at least one, with a
single specific exception. For the case of $\Sigmatil=\matS$ we
also show that for each even degree there are precisely two
exceptions.

\end{abstract}

\section{Problem and new realizability results}\label{new:results:section}
A thorough description of the problem faced below, including
several remarks on its relevance to other areas of mathematics and
a comprehensive survey on the partial solutions obtained over the
time, was given in the first part of the present
paper~\cite{partI}. For this reason we confine ourselves here to
the notation and known facts strictly needed to state and prove
our new results.

\paragraph{Basic definitions}
A \emph{branched covering} is a map $f:\Sigmatil\to\Sigma$, where
$\Sigmatil$ and $\Sigma$ are closed connected surfaces and $f$ is
locally modelled on maps of the form $\matC\ni z\mapsto
z^k\in\matC$ for some $k\geqslant1$. The integer $k$ is called the
\emph{local degree} at the point of $\Sigmatil$ corresponding to
$0$ in the source $\matC$. If $k>1$ then the point of $\Sigma$
corresponding to $0$ in the target $\matC$ is called a
\emph{branching point}. The number of branching points, which is
finite by compactness, will be denoted by $n$. We define $d$ as
the degree of the genuine covering that is obtained by removing
the branching points in $\Sigma$ and all their pre-images in
$\Sigmatil$ and taking the restriction of $f$. We also denote the
number of pre-images of the $i$-th branching point on $\Sigma$ by
$m_i$ and the local degrees at these pre-images by
$(d_{ij})_{j=1}^{m_i}$. It is easy to see that
$(d_{ij})_{j=1}^{m_i}$ forms a partition of $d$. In the sequel we
will always assume that in a partition $(d_1,\ldots,d_m)$ of $d$
we have $d_1\geqslant\ldots\geqslant d_m$.

\paragraph{Branch data}
Suppose we are given closed connected surfaces $\Sigmatil$ and
$\Sigma$, integers $n\geqslant 0$ and $d\geqslant 2$, and for
$i=1,\ldots,n$ a partition $(d_{ij})_{j=1}^{m_i}$ of $d$. The
5-tuple $\big(\Sigmatil,\Sigma,n,d,(d_{ij})\big)$ will be called
the \emph{branch datum} of a candidate branched covering.

\paragraph{Compatibility}
A branch datum is called \emph{compatible} if the following
conditions hold:
\begin{enumerate}
\item\label{RHcond}
$\chi(\Sigmatil)-\ntil=d\cdot(\chi(\Sigma)-n)$;
\item\label{Pcond}
$n\cdot d-\ntil$ is even;
\item\label{OOcond} If $\Sigma$ is
orientable then $\Sigmatil$ is also orientable;
\item\label{NNcond} If $\Sigma$ is non-orientable and $d$ is odd
then $\Sigmatil$ is also non-orientable;
\item\label{ONcond} If
$\Sigma$ is non-orientable but $\Sigmatil$ is orientable then each
partition $(d_{ij})_{j=1}^{m_i}$ of $d$ refines the partition
$(d/2,d/2)$.
\end{enumerate}
The meaning of Condition~\ref{ONcond} is that
$(d_{ij})_{j=1}^{m_i}$ is obtained by juxtaposing two
partitions of $d/2$ and reordering. Note that $d$ is even by
Condition~\ref{NNcond}.

\paragraph{The problem}
It is not too difficult to show that if a branched covering
$\Sigmatil\to\Sigma$ exists then the corresponding branch datum,
with $n,d,(d_{ij}),\ntil$ defined as above, is compatible
(see~\cite{partI} for an explanation of Condition~\ref{Pcond}, the
other ones are obvious). The so-called \emph{Hurwitz existence
problem} is the question of which compatible branch data are
actually realized by some branched covering.

Thanks to the contributions of many authors, among which we will
only mention the fundamental one by Edmonds, Kulkarni and
Stong~\cite{EKS}, we now know that the Hurwitz existence problem
has a positive solution whenever $\chi(\Sigma)\leqslant 0$, and
that the case where $\Sigma$ is the projective plane reduces to
the case where $\Sigma$ is the sphere $\matS$. For this reason we
will assume henceforth that $\Sigma=\matS$.

After being neglected for many years, the problem was picked up
again recently by Bar\'anski~\cite{Baranski}, who proved some
existence results for $\Sigmatil=\Sigma=\matS$, and by
Zheng~\cite{Zheng}, who introduced a new approach to face the
question and obtained many interesting experimental results. In
this paper we extend the technique employed by Bar\'anski to prove
the following results (with $\matT$ denoting the torus in the
second statement and $g\matT$ denoting the orientable surface of
genus $g$ in the third one):

\begin{thm}\label{(d-2,2):sphere:thm}
The compatible and non-realizable branch data of the form
$\big(\matS,\matS,3,d,(d-2,2), (d_{2j}),(d_{3j})\big)$ are
precisely those of the following types:
\begin{itemize}
\item
$\big(\matS,\matS,3,2k,(2k-2,2),(2,\ldots,2),(2,\ldots,2)\big)$
with $k>2$; \item
$\big(\matS,\matS,3,2k,(2k-2,2),(2,\ldots,2),(k+1,1,\ldots,1)\big).$
\end{itemize} In particular, if $d$ is odd then all the data of the
relevant form are realizable.
\end{thm}

\begin{thm}\label{(d-2,2):torus:thm}
With the single exception of
$\big(\matT,\matS,3,6,(4,2),(3,3),(3,3)\big)$, every compatible
branch datum of the form
$\big(\matT,\matS,3,d,(d-2,2),(d_{2j}),(d_{3j})\big)$ is
realizable.
\end{thm}

\begin{thm}\label{(d-2,2):high:genus:thm}
If $g\geqslant 2$ then every compatible branch datum of the form
$\big(g\matT,\matS,3,d,(d-2,2),(d_{2j}),(d_{3j})\big)$
 is realizable.
\end{thm}

\begin{rem}
\emph{The odd-degree part of
Theorem~\ref{(d-2,2):sphere:thm} was already established in an
earlier version of this paper. Later an interesting preprint of
Pakovich~\cite{pako}, using an entirely different approach,
reproved the result and extended it to the case of branch data of
form $\big(\matS,\matS,n,d,(d-k,k),(d_{2j}),\ldots,(d_{nj})\big)$,
which actually contains the whole of our
Theorem~\ref{(d-2,2):sphere:thm}. Given the independence of the
techniques and the fact that they prompt to extensions in
different directions, we have decided to include the proof of
Theorem~\ref{(d-2,2):sphere:thm} anyway.}
\end{rem}

To explain our motivation for considering the case where one of the
partitions has the special form $(d-2,2)$, we mention the fact that
``small'' partitions were already considered in the literature. The
following is for instance known (see~\cite{EKS} and the discussion
in~\cite{partI}):

\begin{thm}\label{full:cycle:thm}
If one of the partitions of $d$ in a compatible branch datum is
given by $(d)$ only then the datum is realizable
\end{thm}

The paper~\cite{EKS} also contains some results for the case where
one of the partitions is $(d-1,1)$ and more generally
$(*,1,\ldots,1)$. In particular, the following statement is
established:

\begin{prop}\label{(d-1,1):prop}
The compatible and non-realizable branch data of the form
$\big(\matS,\matS,n,d,(d-1,1),(d_{2j}),\ldots,(d_{nj})\big)$ are
precisely those of the following types:
\begin{itemize}
\item $\big(\matS,\matS,n,4,(3,1),(2,2),\ldots,(2,2)\big)$ with
$n\geqslant 2$;
\item
$\big(\matS,\matS,3,2k,(2k-1,1),(2,\ldots,2),(2,\ldots,2)\big)$.
\end{itemize}
\end{prop}

Under this perspective, the case $(d-2,2)$ that we consider is
then the next natural one to face.

\section{A geometric criterion for existence}\label{gbara:section}

In this section we prove that the Hurwitz existence problem with
$\Sigma=\matS$ has a geometric equivalent in terms of the
existence of certain families of graphs on the putative covering
surface $\Sigmatil$. Our result extends that proved by
Bar\'anski~\cite{Baranski} for $\Sigmatil=\matS$. A little
machinery has to be developed to give the statement.

\paragraph{Minimal checkerboard graphs}
We begin with a notion also used in~\cite{partI}. We call
\emph{checkerboard graph} a finite $1$-subcomplex of the surface
$\Sigmatil$ whose complement consists of open discs each bearing a
color black or white, so that each edge separates black from
white. We regard checkerboard graphs up to homeomorphism of
$\Sigmatil$ and switching of colors. A checkerboard graph is
called \emph{minimal} if at every vertex all the black germs of
discs are contained in the same global disc, and similarly for
white.

\begin{lemma}
If $G$ is minimal then $\Sigmatil\setminus G$ consists of one
black and one white disc.
\end{lemma}

\begin{proof}
If $x_0$ and $x_1$ are vertices of $G$ joined by an edge, then
there is a black disc incident to both $x_0$ and $x_1$, so there
is only one global black disc incident to $x_0$ and $x_1$, and
similarly for white. Now the conclusion follows from the remark
that $G$ is connected, because its complement consists of
discs.\end{proof}

\paragraph{The case $\Sigmatil=\matS$ and vertices of minimal graphs}
If $\Sigmatil=\matS$ then of course a minimal checkerboard graph
is an embedded circle, as in~\cite{Baranski}, that we consider not
to have any vertex at all. On the contrary, if $\Sigmatil$ has
positive genus, every minimal checkerboard graph has vertices of
valence greater than 2, and we can disregard those of valence 2,
which we will do henceforth.

\paragraph{Minimal graphs from permutations}
We will now prove that for any $g$ there are finitely many minimal
checkerboard graphs on the closed orientable surface $g\matT$ of
genus $g>0$, and that they can be constructed algorithmically. Fix
two copies $B$ and $W$ of the unit disc of $\matC$, and for every
integer $p>1$ denote by $\calR_p^B$ (respectively, $\calR_p^W$)
the set of $p$-th roots of unity on $\partial B$ (respectively,
$\partial W$).

\begin{prop}\label{permu:mingra:prop}
Suppose that $g>0$. Then:
\begin{enumerate}
\item Let $p>1$ and $f\in\permu(\calR_p^B)$. Define
$\ftil\in\permu(\calR_p^B)$ as $\ftil(z)=f(\ee^{2\pi i/p}\cdot
z)$. Suppose that:
\begin{itemize}
\item $f$ consists of $q$ cycles each having length at least $2$;
\item $\ftil$ is a full cycle; \item $p-q=2g$.
\end{itemize}
Then define $h:\partial W\to\partial B$ as
$$h(\ee^{2\pi i(k+t)/p})=\ee^{2\pi i t/p}\cdot (\ftil)^k(1),
\qquad k=0,\ldots,p-1,\quad 0\leqslant t<1,$$ and note that $h$ is
bijective. Define $f'\in\permu(\calR_p^W)$ as $f'=h^{-1}\compo
f\compo h$ and remark that $h$ induces a homeomorphism
$\hbar:\partial B/_{\!f}\to \partial W/_{\!f'}$. Then:
\begin{itemize}
\item The gluing of $B/_{\!f}$ to $W/_{\!f'}$ along $\hbar$ gives
$g\matT$, with $g>0$; \item The image in $g\matT$ of $\partial B$,
which is equal to the image of $\partial W$, is a minimal
checkerboard graph on $g\matT$;
\end{itemize}
\item All minimal checkerboard graphs on $g\matT$ arise
as above for suitable $p>1$ and $f\in\permu(\calR_p^B)$.
\end{enumerate}
\end{prop}

Before turning to the proof of this result, we provide an
alternative description of the statement and an example which show
that the construction, despite its apparent complication, is
actually straight-forward. Let the cycles of $f$ have lengths
$\ell_1,\ldots,\ell_q$. Give labels $a^{(i)}_j$ to the points of
$\calR_p^B$ for $i=1,\ldots,q$ and $j=0,\ldots,\ell_i-1$, in such
a way that $f(a^{(i)}_j)=a^{(i)}_{j+1}$ with indices modulo
$\ell_i$. An example is shown for
$p=7,q=3,\ell_1=\ell_2=2,\ell_3=3$ in
Fig.~\ref{mingra-example1:fig}-left,
    \begin{figure}
    \begin{center}
    \input{new_mingra-example1.pstex_t} 
    \mycap{A permutation of the black $7$-th roots of unity and the associated white permutation.} \label{mingra-example1:fig}
    \end{center}
    \end{figure}
where for simplicity we have replaced $a^{(1)},a^{(2)},a^{(3)}$ by
$a,b,c$. The definition of $\ftil$ is then as follows: from any
point of $\calR_p^B$, first go to the next point on $\partial B$
in the counter-clockwise order, and then jump to the next point
according to $f$. In our case, we get the full cycle $a_0\mapsto
b_0\mapsto a_1\mapsto b_1\mapsto c_2\mapsto c_1\mapsto c_0$. If we
label the points of $\calR_p^W$ according to the labels in this
cycle $\ftil$, as in Fig.~\ref{mingra-example1:fig}-right, the
definitions of $f'$ and $h$ are now immediately understood: again
$f'(a^{(i)}_j)=a^{(i)}_{j+1}$, while $h:\partial W\to\partial B$
maps each open arc of $\partial W\setminus\calR_p^W$ to the open
arc of $\partial B\setminus\calR_p^B$ having \emph{first} endpoint
with the same label.  See Fig.~\ref{mingra-example1bis:fig} for
the identifications of arcs in our example, and
    \begin{figure}
    \begin{center}
    \input{new_mingra-example1bis.pstex_t} 
    \mycap{Example continued: identifications of arcs on $\partial B$ and $\partial W$.} \label{mingra-example1bis:fig}
    \end{center}
    \end{figure}
Fig.~\ref{mingra-example1ter:fig} for the resulting checkerboard
coloring of the genus-2 surface.
    \begin{figure}
    \begin{center}
    \input{new_mingra-example1ter.pstex_t} 
    \mycap{Example continued: a minimal checkerboard coloring of the genus-2 surface.} \label{mingra-example1ter:fig}
    \end{center}
    \end{figure}

\begin{proof}
Suppose a minimal checkerboard graph $G$ is given. Since $g$ is
positive, $G$ has vertices. Suppose there are $q$ of them, and
label them by $a^{(i)}$ for $i=1,\ldots,q$. If $a^{(i)}$ has
valence $2\ell_i$, label the germs of white and black discs around
$a^{(i)}$ in positive order as
$$a^{(i)}_0,a^{(i)}_0,a^{(i)}_1,a^{(i)}_1,\ldots,
a^{(i)}_{\ell_i-1},a^{(i)}_{\ell_i-1}$$ starting from an arbitrary
\emph{white} germ. Recalling that there is only one white disc,
label the edges of $G$ as $\alpha_1,\ldots,\alpha_p$ while
following the boundary of the white disc in the \emph{negative}
fashion. Now read the labels of the germs of discs and of the
edges on the boundary of the abstract white and black discs $W$
and $B$ with respect to the obvious immersions from $W$ and $B$ to
the surface. The maps $f$ and $f'$ of the definition are now those
telling that the vertices $a^{(i)}_j$, for $j=0,\ldots,\ell_i-1$,
become the single vertex $a^{(i)}$ in the surface, and describing
the way the germs are cyclically arranged around $a^{(i)}$.
Similarly, the map $h$ says that each arc $\alpha_k$ on $\partial
W$ is glued to the corresponding $\alpha_k$ on $\partial B$. The
cell decomposition of the surface associated to $G$ consists of
$q$ vertices, $p$ edges, and 2 discs, whence $p-q=2g$, which
implies that $G$ arises from $f\in\permu(\calR^B_p)$ as described
in point 1.

Now suppose $f\in\permu(\calR^B_p)$ is given. For each orbit
$a^{(i)}_0,\ldots,a^{(i)}_{\ell_i-1}$ of $f$ we identify the
points of the orbit to a single point $a^{(i)}$ and we give to the
quotient space a planar structure near $a^{(i)}$ so that the germs
of discs arising from the $a^{(i)}_j$'s are cyclically arranged in
the order $j=0,\ldots,\ell_i-1$. With this choice a system of
attaching loops for white discs is well-defined. The condition
that $\ftil$ is a full cycle translates the fact that a single
white disc is glued to the black disc, and the condition that
$p-q=2g$ means that the result of the gluing is the orientable
genus-$g$ surface.
\end{proof}

\begin{cor}\label{finite:mingra:cor}
On $g\matT$ there are finitely many minimal checkerboard graphs.
\end{cor}

\begin{proof}
With the notation of the previous proposition, we have that
$q\leqslant p/2$ from the assumption that $f$ has no fixed points.
Whence $p\leqslant 4g$ from the condition $p-q=2g$. The conclusion
follows.
\end{proof}

\paragraph{Genus $0$ and $1$}
As already noticed, the only minimal checkerboard graph on $\matS$ is a
plain circle, which splits $\matS$ into an embedded black and an
embedded white disc. This
corresponds in Proposition~\ref{permu:mingra:prop} to the limit
case $p=q=0$.

For genus $1$, the argument proving
Corollary~\ref{finite:mingra:cor} shows that $p\leqslant 4$.
Moreover $p-q=2$. Now we note the following easy general fact:

\begin{rem}\emph{Let $f\in\permu(\calR_p^B)$ and suppose that
$f$ maps some point in $\calR_p^B$ to the point preceding it in
the counter-clockwise order on $\partial B$. Then the
corresponding $\ftil$ as in Proposition~\ref{permu:mingra:prop}
has a fixed point, so it is not a full cycle}.\end{rem}

\begin{figure}
    \begin{center}
    \includegraphics[scale=0.5]{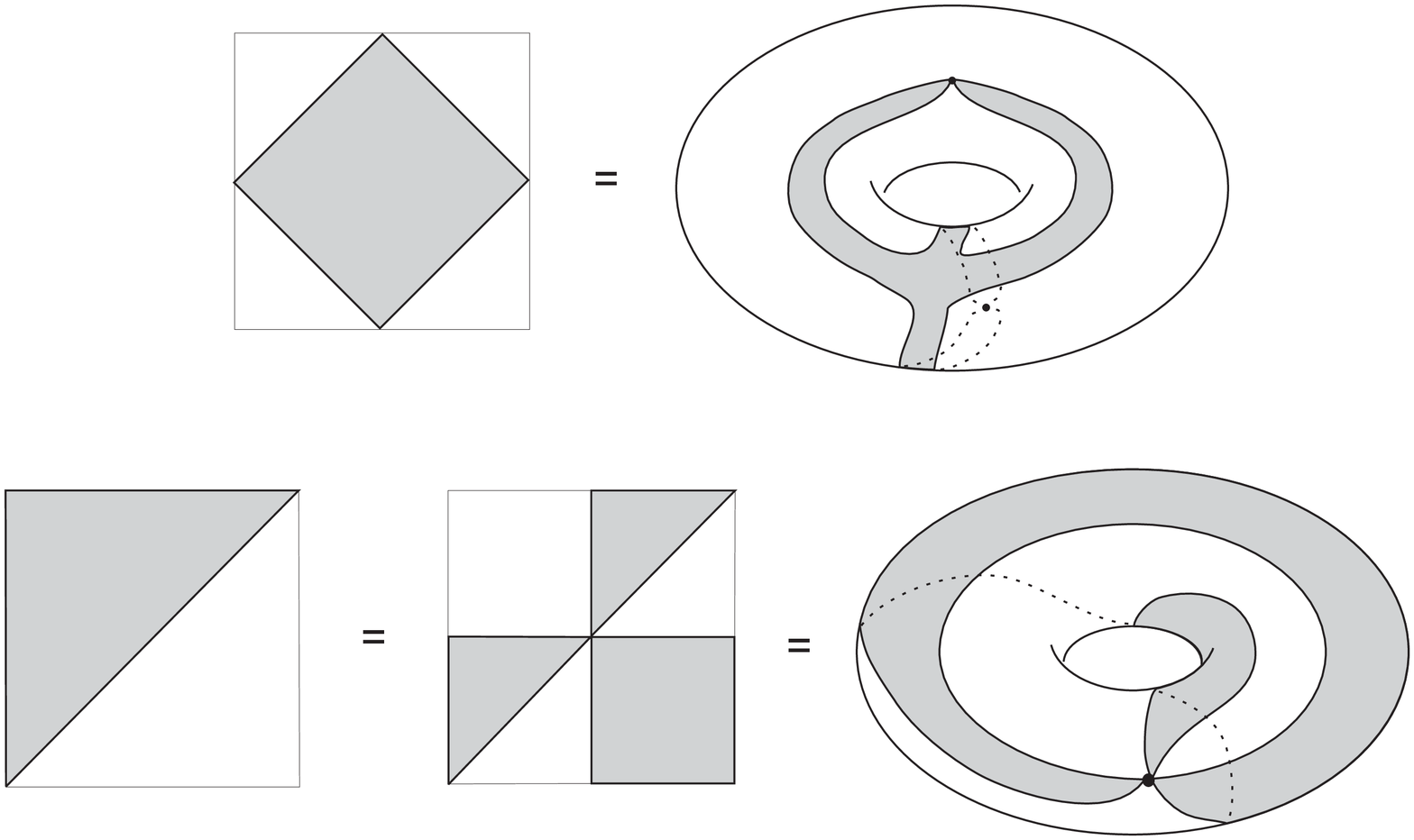}
    \mycap{Minimal checkerboard graphs on the torus.} \label{mingra-torus:fig}
    \end{center}
    \end{figure}

According to this result, only the $p$'s and $f$'s as in
Fig.~\ref{permu-for-torus:fig}
    \begin{figure}
    \begin{center}
    \input{new_permu-for-torus.pstex_t} 
    \mycap{Permutations of roots of unity relevant for the torus.} \label{permu-for-torus:fig}
    \end{center}
    \end{figure}
can be relevant for $g=1$. Both these permutations $f$ satisfy the
condition that the corresponding $\ftil$ as in
Proposition~\ref{permu:mingra:prop} is a full cycle. The
associated minimal checkerboard graphs (and colorings) of the
torus are shown in Fig.~\ref{mingra-torus:fig}.

\paragraph{Existence of coverings} We now prove our extension of
Bar\'anski's criterion for existence of a branched covering onto
$\matS$. We denote by $V(G)$ the set of vertices of a graph $G$
and by ${\rm val}_G(w)$ the valence in $G$ of a vertex $w$.

\begin{thm}\label{gTtoS:thm}
A branch datum $\big(g\matT,\matS,n,d,(d_{ij})\big)$ is
realizable if and only if there exists on $g\matT$ a minimal
checkerboard graph $G$ associated to $f\in\permu(\calR^B_p)$ as in
Proposition~\ref{permu:mingra:prop}, a set $S$ of $n\cdot d-p$
points of $G\setminus V(G)$, a labelling in $\{1,\ldots,n\}$ of
$V(G)\cup S$ and a collection $\Gamma=\bigcup\Gamma_{ij}$ of
graphs on $g\matT$ such that:
\begin{itemize}
\item The $\Gamma_{ij}$'s are pairwise disjoint trees; \item The
set of vertices of $\Gamma$ is equal to $G\cap\Gamma$ and is also
equal to $V(G)\cup S$; \item Pulling back the labelling of the
points through the projection $\partial B\to G$ (or, equivalently,
$\partial W\to G$) we find $1,\ldots,n,\ldots,1,\ldots,n$, with
each string $1,\ldots,n$ appearing $d$ times, in counter-clockwise
order; \item $\Gamma_{ij}\cap G$ consists of points having label
$i$; \item The number of edges of $\Gamma_{ij}$ is given by
$$d_{ij}-1-\sum\limits_{w\in V(G)\cap\Gamma_{ij}}\left(\frac12{\rm
val}_G(w)-1\right).$$
\end{itemize}
\end{thm}

\begin{proof}
Suppose the labelling of points and the trees $\Gamma_{ij}$ exist.
Let us assign to all the points of $\Gamma_{ij}$ the label $i$,
and let us call $i$-\emph{portion} of a set $X$ a connected
component of the set of points of $X$ having label $i$. We claim
the following assertion:
\begin{enumerate}
\item[] \emph{The complement in $g\matT$ of $G\cup\Gamma$ consists
of $d$ black and $d$ white discs. On each black (respectively,
white) disc the labelled portions of $\Gamma$ cyclically found on
the boundary in the positive (respectively, negative) order have
labels $1,\ldots,n$.}
\end{enumerate}

Assuming the assertion for a moment, let us explain how to
construct the desired covering.  Since the $\Gamma_{ij}$'s are
disjoint trees, we can collapse them to points without changing
the topology of $g\matT$. After this collapse, $G$ gives a graph
$G'$ on $g\matT$ with labelled vertices and the complement
consisting of $d$ black and $d$ white discs. Moreover on each
black (respectively, white) disc the labels of the vertices
cyclically found in the positive (respectively, negative) order on
the boundary are $1,\ldots,n$. And each edge of $G'$ separates
discs of different colors, because this is true for $G$. The last
two statements easily imply that the discs of $g\matT\setminus G'$
have embedded closures.

Consider now the analogous decomposition of $\matS$ as shown in
Fig.~\ref{sphere-deco:fig}, with one black and one white disc. We
map the closure of each black
    \begin{figure}
    \begin{center}
    \input{new_sphere-deco.pstex_t} 
    \mycap{The sphere as a union of a black and a white disc.} \label{sphere-deco:fig}
    \end{center}
    \end{figure}
(respectively, white) disc on $g\matT\setminus G'$ to the closure
of the black (respectively, white) disc on $\matS$, matching the
labels of the vertices. These maps can be arranged to coincide on
the edges of $G$ (which are shared by white and black discs). The
result is of course a branched covering, and to conclude we only
need to prove that the local degree at the point of $g\matT$
arising from the collapse of $\Gamma_{ij}$ is $d_{ij}$. This local
degree is half of the valence of the vertex in $G'$. We now show
that half of this valence is given by the number of edges of
$\Gamma_{ij}$ plus
$$1+\sum\limits_{w\in V(G)\cap\Gamma_{ij}}\left(\frac12{\rm
val}_G(w)-1\right)$$ which implies the conclusion by the
assumption. The formula is of course correct when $\Gamma_{ij}$
has no edges, because it must consist of a single point of $S\cup
V(G)$. Adding an edge with endpoint $w$ to an already existing
$\Gamma_{ij}$, half of the valence in $G'$ increases by
$\frac12{\rm val}_G(w)=1+\left(\frac12{\rm val}_G(w)-1\right)$.
The conclusion follows.

Let us now prove the claimed assertion. Note first that, by the
second assumption and the last one, the total number of edges in
$\Gamma$ is
$$\sum_{i=1}^n\sum_{j=1}^{m_i}(d_{ij}-1)-\sum_{w\in
V(G)}\left(\frac12{\rm val}_G(w)-1\right).$$ The first sum gives
$n\cdot d-\ntil$ while the second sum gives $p-q$. Since $n\cdot
d-\ntil=2d-2+2g$ and $p-q=2g$, the total number of edges is
$2d-2$. Before adding edges, \emph{i.e.} in $g\matT\setminus G$,
there are two discs, and each time we add an edge we increase by
one the number of discs, therefore the total number of discs in
$g\matT\setminus G'$ is indeed $2d$.

To proceed, let us call $G$-\emph{arc} a connected component of
$G\setminus (S\cup V(G))$ and let us assign the label $i$ to all
the points of a $G$-arc having ends $i$ and $i+1$ in
counter-clockwise order with respect to $B$. If $U$ is a component
of $g\matT\setminus(G\cup\Gamma)$ and we are given distinct points
of $\partial U$, we can speak of the two halves into which the
points split $\partial U$, even if we do not know that the closure
of $U$ is embedded. We claim the following: \emph{for every
component $U$ of $g\matT\setminus (G\cup\Gamma)$, for all $i\neq
i'$ and for every two distinct points of label $i'$ in $\partial
U$, at least one of the halves into which the points split
$\partial U$ contains a $G$-arc of label $i$. The other half is
either entirely labelled with $i'$ or it also contains a $G$-arc
of label $i$}. We prove this claim again by imagining that
$\Gamma$ has been constructed adding one edge after the other. The
claim is of course true before adding edges, \emph{i.e.} in $G$.
Assume the claim is true up to some stage. When we add the next
edge inside a certain region $U$, we replace $U$ by two regions
$U'$ and $U''$, and both $\partial U'$ and $\partial U''$ differ
from $\partial U$ as follows: a (possibly immersed) segment having
ends of a certain label $i'$ is replaced by an embedded segment
entirely labelled by $i'$. One easily sees that the property
stated in the claim is preserved by such a transformation (recall
that $\Gamma$ consists of trees).

The claim implies that each component of $g\matT\setminus
(G\cup\Gamma)$ contains at least one $G$-arc of each label. But
there are $2d$ components, $n$ labels and $n\cdot d$ arcs, and
each $G$-arc is shared by 2 components, so the boundary
of each component of $g\matT\setminus
(G\cup\Gamma)$ contains precisely one $G$-arc for each label. For
the same reasons, there are $d$ black and $d$ white components.
The cyclic order of the labels of the $G$-arcs is preserved at
each insertion of an edge, so it is $1,\ldots,n$ in positive
(respectively, negative) order on the black (respectively, white)
discs. To conclude the proof of the assertion we only need to note
that a $G$-arc of label $i$ and one of label $i+1$ are always
separated by an $i$-portion of $\Gamma$. The assertion, and hence
the existence of the covering, are eventually established.

Having proved the sufficiency of the conditions stated in the
theorem, let us turn to the opposite implication. We then assume
that the desired branched covering $g\matT\to\matS$ exists and we
suppose, which of course we can, that the branching points on
$\matS$ are the labelled points in Fig.~\ref{sphere-deco:fig}.

Let $H$ be the pre-image in $g\matT$ of the circle in
Fig.~\ref{sphere-deco:fig} and $\Delta=(\Delta_{ij})$ be the
collection of pre-images of the branching points. We give a label
$i$ to $\Delta_{ij}$ and a black or white color to each disc in
$g\matT\setminus H$, according to the color of its projection. We
consider moves on the pair $(H,(\Delta_{ij}))$ which preserve the
following properties:
\begin{itemize}
\item $H$ is a checkerboard graph on $g\matT$; \item The
$\Delta_{ij}$'s are mutually disjoint trees; \item The set of
vertices of each $\Delta_{ij}$ is given by its intersection with
$H$; \item The number of edges of $\Delta_{ij}$ plus a
contribution of $\frac12{\rm val}_H(v)-1$ for each vertex $v$ of
$H$ contained in $\Delta_{ij}$ is $d_{ij}-1$.
\end{itemize}
Note that the conditions are of course met at the beginning. There
are two moves that we employ, a black one and a white one. The
white one applies if at some vertex of $H$ there are two distinct
white discs.  The move itself is described in
Fig.~\ref{white-move:fig}
    \begin{figure}
    \begin{center}
    \input{new_white-move.pstex_t} 
    \mycap{The white move.} \label{white-move:fig}
    \end{center}
    \end{figure}
where we assume that $R_1\neq R_2$. Note that there is some
arbitrariness in the move, related to the choice of the vertex
where it is applied, to the white discs of $g\matT\setminus H$
which are merged, and to the new position of the edges in $\Delta$
coming from the previous steps. The black move is completely
analogous to the white one.

Let us now apply the moves as long as possible and call
$(G,(\Gamma_{ij}))$ the final $(H,(\Delta_{ij}))$. By definition
$G$ is a minimal checkerboard graph (to ensure the minimality, we
stipulate that the vertex set of $G$ consists of the points with
link of cardinality at least three). Let us define $S$ to be the
set of points of $G\cap\Gamma$ which are not vertices of $G$, and
let us assign the label $i$ to a point of $S\cup V(G)$ if it
belongs to some $\Gamma_{ij}$. It is now a routine matter to
repeat the above arguments in the reverse order to show that the
assumptions of the theorem are verified.\end{proof}

\paragraph{Useful minimal checkerboard graphs}
According to Theorem~\ref{gTtoS:thm}, only the minimal
checkerboard graphs supporting the further structures described in
the theorem are relevant to the Hurwitz existence problem, so
\emph{a priori} some graph may be dispensable. However we have the
following:

\begin{prop}\label{all:mingra:useful:prop}
If $G$ is a minimal checkerboard graph for $g\matT$ then there
exists a branch datum with $n=3$ and suitably large $d$ which is
realized according to Theorem~\ref{gTtoS:thm} using the graph $G$.
\end{prop}

\begin{proof}
Let us give label $1$ to all the vertices of $G$. Let $S$ consist
of two points (with labels $2$ and $3$) on each edge of $G$. Now
consider in $B$ and $W$ the trees as shown in
Fig.~\ref{all-mingra-useful:fig}.
    \begin{figure}
    \begin{center}
    \input{new_all-mingra-useful.pstex_t} 
    \mycap{Trees showing that all minimal checkerboard graphs give rise
    to some covering.} \label{all-mingra-useful:fig}
    \end{center}
    \end{figure}
These trees give a branched covering of type $(d),(d),(*)$.
\end{proof}

This result does not imply, however, that all minimal checkerboard
graphs are necessary to construct the realizable branch data,
because a branch datum could be constructed using different
graphs. The next remark shows that this is what happens in the
case of the torus.

\begin{rem}
\emph{Let us denote by $G_1$ and $G_2$ the minimal checkerboard
graphs on the torus that are shown in
Fig.~\ref{mingra-torus:fig}-top and in
Fig.~\ref{mingra-torus:fig}-bottom, respectively. The graph $G_1$,
labels and trees in Fig.~\ref{torus-useful:fig} give a covering
    \begin{figure}
    \begin{center}
   \input{revised_new_torus-useful.pstex_t} 
    \mycap{Realization of a covering $\matT\to\matS$.} \label{torus-useful:fig}
    \end{center}
    \end{figure}
realizing the datum $\big(\matT,\matS,4,2,(2),(2),(2),(2)\big)$,
whose non-trivial automorphism is the hyperelliptic involution,
which cannot be obtained using the graph $G_2$. On the contrary,
we will now show that any datum that can be realized using $G_2$
can also be realized using $G_1$. Indeed, suppose that a
realization of some datum is given by a family of graphs with
labelled vertices $G_2\cup\Gamma\subset\matT$, according to
Theorem~\ref{gTtoS:thm}. We first claim that in the white disc $W$
there is an edge $e$ with ends in internal points of different
edges of $G_2$. If this is not the case then there is a disc in
$W\setminus(G_2\cup\Gamma)$ doubly incident to $v$, but we have
shown within the proof of Theorem~\ref{gTtoS:thm} that (even after
collapsing each component of $\Gamma$ to a point) these discs have
embedded closures. The claim is proved and we can now collapse $e$
to a point and perform a white move at $v$, which turns $G_2$ into
$G_1$ and $\Gamma$ into another union of trees realizing the same
covering.}
\end{rem}

\section{Dessins d'enfants}\label{dessins:section}
In this short section we review a technique introduced
in~\cite{partI} to investigate existence of branched coverings. It
is based on the notion of dessin d'enfant due to
Grothendieck~\cite{Groth}. This notion in its original form was
only relevant to the case of $n=3$ branching points, but
in~\cite{partI} it was extended to arbitrary $n\geqslant 4$. In
this paper, however, we will only need it for $n=3$, so we quickly
present it in its simplified form here.

To begin we recall that a graph (a finite 1-complex) is called
\emph{bipartite} if its set of vertices is split as $V_1\sqcup
V_2$ and each edge has one endpoint
in $V_1$ and one in $V_2$.

\begin{defn}
\emph{A \emph{dessin d'enfant} on the surface $\Sigmatil$ is a
bipartite graph $D\subset\Sigmatil$ such that $\Sigmatil\setminus
D$ consists of open discs. The \emph{length} of one of these discs
is the number of edges of $D$ along which its boundary passes
(with multiplicity).}
\end{defn}

\begin{prop}\label{dessins:coverings:prop}
The realizations of a branch datum
$\big(\Sigmatil,\matS,3,d,(d_{ij})\big)$ correspond to the dessins
d'enfants $D\subset\Sigmatil$ with set of vertices split as
$V_1\sqcup V_2$ such that for $i=1,2$ the vertices in $V_i$ have
valences $(d_{ij})_{j=1}^{m_i}$, and the discs in
$\Sigmatil\setminus D$ have lengths $(2d_{3j})_{j=1}^{m_3}$.
\end{prop}

\begin{proof}
Suppose a realization $f:\Sigmatil\to\matS$ exists, let the
branching points be $p_1,p_2,p_3$, choose an arc $\alpha$ joining
$p_1$ to $p_2$ and avoiding $p_3$, and define $D$ as
$f^{-1}(\alpha)$.  Setting $V_i=f^{-1}(p_i)$ for $i=1,2$, it is
clear that $D$ is a bipartite graph with partition of vertices
$V_1\sqcup V_2$ with the required valences. Now
$\matS\setminus\alpha$ is an open disc and the restriction of $f$
to any component of $\Sigmatil\setminus D$ is a covering onto this
disc with a single branching point.  Such a covering is always
modelled on the covering $z\mapsto z^k$ of the open unit disc onto
itself, so the components of $\Sigmatil\setminus D$ are open
discs. More precisely, there is one such disc for each element of
$f^{-1}(p_3)$, and it is easy to see that the $j$-th one has
length $2d_{3j}$.

Reversing this construction is a routine matter that we can leave
to the reader.\end{proof}

To apply Proposition~\ref{dessins:coverings:prop} we will
often switch the viewpoint: instead
of trying to embed a dessin $D$ in the surface $\Sigmatil$, we
will try to thicken a given bipartite graph $D$ with certain
prescribed valences of vertices to a surface with boundary, so to
get $\Sigmatil$ by capping off the boundary circles.

\section{Existence results}\label{exist:section}

In this section we establish Theorems~\ref{(d-2,2):sphere:thm}
to~\ref{(d-2,2):high:genus:thm}. Our proofs will repeatedly employ
Theorem~\ref{full:cycle:thm}. We begin with some technical
notions.

\paragraph{Diagrams and accessible trees} To investigate the
branch data involving the partition $(d-2,2)$ we will apply the
geometric existence criterion Theorem~\ref{gTtoS:thm}, of which we
will use the notation throughout. For the sake of brevity, a
triple $G,S,\Gamma$ as in the statement of this theorem will be
called a \emph{diagram} realizing the corresponding branch datum.
A connected component of $\Sigmatil\setminus(G\cup\Gamma)$ will be
called a \emph{face} of this diagram. Recalling that each tree
$\Gamma_{ij}\in\Gamma$ corresponds to some entry $d_{ij}$ in one
of the partitions of the relevant branch datum, we will call
\emph{degree} of $\Gamma_{ij}$ the associated $d_{ij}$. The
collection of degrees of the trees incident to a given face
(\emph{i.e.} having non-empty intersection with its closure in
$\Sigmatil$) will be called the \emph{type} of the face.

For the rest of the paper we fix branch data as in the assumptions
of Theorems~\ref{(d-2,2):sphere:thm}, \ref{(d-2,2):torus:thm},
and~\ref{(d-2,2):high:genus:thm}, \emph{i.e.} we will henceforth
assume that $m_1=2$ and $(d_{11},d_{12})=(d-2,2)$. Moreover, in
all our diagrams the trees $\Gamma_{11}$ and $\Gamma_{12}$ will
correspond to this partition of $d$. For $i=2,3$ we will say that
a tree $\Gamma_{ij}$ of a diagram is \emph{accessible} if there is
a $G$-arc with one endpoint in $\Gamma_{ij}$ and the other
endpoint in $\Gamma_{11}$. Recall that a $G$-arc is a connected
component of $G\setminus(S\cup V(G))$. A diagram is called
\emph{accessible} if all $\Gamma_{ij}$'s with $i=2,3$ are
accessible.

\begin{lemma}\label{accessib:diag:lem}
In a diagram realizing our branch datum there is at most one
inaccessible tree, which must be a point that is not a vertex of $G$.
\end{lemma}

\begin{proof}
Recall the collapse used in the proof of Theorem~\ref{gTtoS:thm},
where each $\Gamma_{ij}$ is shrunk to a point $v_{ij}$ labelled
$i\in\{1,2,3\}$, and $G$ is transformed to some graph $G'$. The
valence of $v_{ij}$ is $2d_{ij}$, so $v_{11}$ has valence $2d-4$
and $v_{12}$ has valence $4$. The complement of $G'$ consists of
$2d$ triangles, each with embedded closure and vertices labelled
$1,2,3$. So there are $2d-4$ triangles incident to $v_{11}$, and,
if $\Gamma_{ij}$ is inaccessible, then $v_{ij}$ does not belong to
any of them. So it belongs to the interior of the union $U$ of the
other $4$ triangles, which are all incident to $v_{12}$. The only
two possibilities for $U$ are as shown in
Fig.~\ref{inaccess_graph:fig},
    \begin{figure}
    \begin{center}
\input{revised_new_inaccess_graph.pstex_t} 
    \mycap{The triangles incident to $v_{12}$.} \label{inaccess_graph:fig}
    \end{center}
    \end{figure}
and the conclusion follows.\end{proof}

\paragraph{Odd-degree coverings by the sphere} In the case
$\Sigmatil=\matS$ it will be convenient to investigate the
data with odd degree $d$ first. Namely, we establish the following:

\begin{prop}\label{(d-2,2):sphere:odd:prop}
If $d$ is odd then every compatible branch datum of the form
$\big(\matS,\matS,3,d,(d-2,2),(d_{2j}),(d_{3j})\big)$ is
realizable.
\end{prop}

\begin{proof}  Note first that compatibility means that $m_2+m_3=d$
and recall that, in a diagram $G,S,\Gamma$ realizing the datum,
$G$ is a plain circle, $S$ consists of $3d$ points, and
$\Gamma_{ij}$ is a tree with $d_{ij}-1$ edges. So $\Gamma_{11}$
has $d-3$ edges and $\Gamma_{12}$ is just one edge.

By Theorem~\ref{full:cycle:thm} the datum is realizable if
$(d_{ij})_{j=1}^{m_i}=(d)$ for some $i\in\{2,3\}$. We will
then assume henceforth that this is not the case, which easily
implies that $d\geqslant 5$ (recall that $d$ is odd). To prove the
proposition the general strategy will
now be to proceed by induction on $d$, starting from the diagrams
$D_1,\ldots,D_4$ of Fig.~\ref{degree5:fig},
    \begin{figure}
    \begin{center}
    \input{new_degree5.pstex_t} 
    \mycap{The diagrams $D_1,\ldots,D_4$, with $\Gamma_{11}$ and $\Gamma_{12}$
    marked by thicker lines. In $D_1,D_2,D_3$ a small arrow
    indicates the only inaccessible vertex, while the $G$-arc $e$
    on $D_4$ will be used later in the proof. The partitions
    $(d_{ij})_{i=2,3}^{j=1,\ldots,m_i}$ are $(4,1),(3,1,1)$ for
    $D_1$, $(4,1),(2,2,1)$ for $D_2$, $(3,2),(3,1,1)$ for $D_3$,
    and $(3,2),(2,2,1)$ for $D_4$.} \label{degree5:fig}
    \end{center}
    \end{figure}
which realize all the relevant data for $d=5$, and successively
applying certain moves to these diagrams so to realize all the
relevant data for arbitrary $d$. However, we will have to deal
with special sorts of data by separate inductions.

We will employ three moves $\mu_1,\mu_2,\mu_3$, that we now
describe:

\begin{itemize}
\item[$\mu_1$:] Choose $G$-arcs $e_2$ and $e_3$ such that each
$e_i$ has one end in $\Gamma_{11}$ and the other end in some
$\Gamma_{ij_i}$, and perform the modifications of
Fig.~\ref{move1:fig};
    \begin{figure}
    \begin{center}
    \input{new_move1_on_diag.pstex_t} 
    \mycap{The move $\mu_1$.}\label{move1:fig}
    \end{center}
    \end{figure}
\item[$\mu_2$:] Choose a $G$-arc $e$ with one end in $\Gamma_{11}$
and one in some $\Gamma_{ij_i}$, and perform the modification of
Fig.~\ref{move2:fig};
    \begin{figure}
    \begin{center}
    \input{new_move2_on_diag.pstex_t} 
    \mycap{The move $\mu_2$.}\label{move2:fig}
    \end{center}
    \end{figure}
\item[$\mu_3$:] Choose a $G$-arc $e$ with one end in $\Gamma_{11}$
and one in some $\Gamma_{ij_i}$, and perform the modification of
Fig.~\ref{move3:fig}.
    \begin{figure}
    \begin{center}
    \input{new_move3_on_diag.pstex_t} 
    \mycap{The move $\mu_3$.}\label{move3:fig}
    \end{center}
    \end{figure}
Since we will sometimes use a systematic iteration of this move, a
new $G$-arc is specified and also labelled $e$ in the new diagram.
\end{itemize}

It is evident that all three moves transform a diagram of degree
$d$ to a diagram of degree $d+2$ which still satisfies the
conditions of Theorem~\ref{gTtoS:thm}. We call a diagram {\em
constructible} if it is obtained from one of the diagrams
$D_1,\ldots,D_4$ by successive application of moves $\mu_*$. We
will prove that each relevant branch datum is realized by a
constructible diagram.

Let us note now that the effect of the $\mu_*$'s on the partitions
$(d_{2j})$ and $(d_{3j})$ is as follows:

\begin{itemize}

\item[$\widehat{\mu}_1$:] Choose $j_2,j_3$, replace $d_{2j_2}$ and
$d_{3j_3}$ by $d_{2j_2}+1$ and $d_{3j_3}+1$ respectively, add a
$1$ at the end of both partitions, and reorder (if necessary:
recall that our partitions are arranged in non-increasing order);

\item[$\widehat{\mu}_2$:] Choose $i\in\{2,3\}$ and $j_i$, let
$\{2,3\}=\{i,i'\}$, replace $d_{ij_i}$ by $d_{ij_i}+2$ and
reorder, and add two $1$'s at the end of $(d_{i'j})$;

\item[$\widehat{\mu}_3$:] Choose $i\in\{2,3\}$ and $j_i$, let
$\{2,3\}=\{i,i'\}$, replace $d_{ij_i}$ by $d_{ij_i}+1$, add a $1$
at the end, add a $2$ at the end of $(d_{i'j})$, and reorder.

\end{itemize}
Note also that the $\widehat{\mu}_*$'s can be viewed as moves
transforming a relevant branch datum of degree $d$ into one of
degree $d+2$. However, \emph{a priori} not every branch move
$\widehat{\mu}_*$ on the branch datum realized by a constructible
diagram is induced by the corresponding geometric move $\mu_*$ on
that diagram. When this happens, we will say that the move
$\widehat{\mu}_*$ itself is \emph{constructible}.  Of course a
move $\widehat{\mu}_*$ applied to indices $j_2,j_3$ (for
$\widehat{\mu}_1$) or $j_i$ (for $\widehat{\mu}_2$ or
$\widehat{\mu}_3$) is constructible if and only if the
corresponding trees $\Gamma_{2j_2}$ and $\Gamma_{3j_3}$ or
$\Gamma_{ij_i}$ are accessible. Lemma~\ref{accessib:diag:lem} then
implies that \emph{a branch move $\widehat{\mu}_*$ is
constructible if applied to indices such that the corresponding
 $d_{2j_2}$ and $d_{3j_3}$ or
$d_{ij_i}$ are larger than $1$.} We also note that a move $\mu_*$
never creates inaccessible trees. Therefore \emph{a branch move
$\widehat{\mu}_*$ is constructible if it is applied to a diagram
constructed starting from $D_4$ or to tree(s) created by previous
geometric moves.} These facts will be used repeatedly below.

We now state and prove a series of claims that will lead to the
conclusion. We start with an easy arithmetic one.

\medskip

\noindent\textsc{Claim 1.} \emph{Either
$(d_{ij})_{i=2,3}^{j=1,\ldots,m_i}$ includes two $1$'s or
$(d_{2j})_{j=1}^{m_2}=(3,2,\ldots,2)$ and
$(d_{3j})_{j=1}^{m_3}=(2,\ldots,2,1)$, up to permutation.}
Assume first there are no $1$'s. Since $d$ is odd we deduce that
$d_{21},d_{31}\geqslant 3$. Therefore
$$\Big(d\geqslant 3+2(m_i-1) \Rightarrow
m_i\leqslant\frac{d-1}{2},\ i=2,3\Big)\Rightarrow m_2+m_3\leqslant
d-1,$$ which gives a contradiction. Suppose then that $d_{3m_3}=1$
and there are no other $1$'s. Then again
$m_2\leqslant\frac{d-1}{2}$, and $m_3\leqslant\frac{d+1}{2}$ by a
similar argument. So both inequalities are equalities and the
conclusion readily follows.

\medskip

We are now ready to proceed with the main part of the proof, where
we show that every relevant branch datum, \emph{i.e.} one of the
form
\begin{equation}\label{S:S:odd:d}
\big(\matS,\matS,3,d,(d-2,2),(d_{2j})_{j=1}^{m_2},(d_{3j})_{j=1}^{m_3}\big)
\end{equation}
with odd $d\geqslant 7$ and $m_2,m_3\geqslant 2$, is realized by
some constructible diagram. The proof is separate for branch data
that involve a partition of the form $(3,\ldots,3,2,\ldots,2)$,
\emph{i.e.} those of form
\begin{equation}\label{S:S:odd:(3,2):d}
\big(\matS,\matS,3,d,(d-2,2),(\underbrace{3,\ldots,3}_k,2,\ldots,2),
(d_{3j})_{j=1}^{m_3}\big).
\end{equation}

\medskip

\noindent\textsc{Claim 2.} \emph{If
$(d_{2j})=(3,\ldots,3,2\ldots,2)$ and
$(d_{3j})=(2,\ldots,2,1,\ldots,1)$ then the datum can be realized
by an accessible constructible diagram.} Since the total number of
entries in the two partitions is $d$, the number $k$ of $3$'s in
the first partition is equal to the number of $1$'s in the second
one, and $k$ is odd. We prove the claim by induction on $k$. The
diagram obtained from $D_4$ by $(d-5)/2$ successive applications
of the move $\mu_3$ starting with the $G$-arc $e$ indicated in
Fig.~\ref{degree5:fig} gives the base $k=1$ of the induction.

For the inductive step, notice that to get the
datum~(\ref{S:S:odd:(3,2):d}) at level $k$ in degree $d$ with the
relevant $(d_{3j})$ we can apply a move $\widehat{\mu}_1$ (with
$j_2=k-1$ and $j_3=m_3-3$) followed by a move $\widehat{\mu}_2$
(with $i=2$ and $j_2=m_2$) to the same datum at level $k-2$ in
degree $d-4$, and the inductive assumption easily implies the
conclusion.

\medskip

\noindent\textsc{Claim 3.} \emph{If $(d_{2j})=(3,2\ldots,2)$ then
the datum can be realized by a constructible diagram.} By
induction on $d\geqslant 5$. For $d=5$ the conclusion is given by
diagrams $D_3$ and $D_4$. For the inductive step, suppose first
that $d_{31}\leqslant 2$. It follows that
$(d_{3j})=(2,\ldots,2,1)$, and this case was settled in Claim 2.
Then we assume $d_{31}\geqslant 3$. So the
datum~(\ref{S:S:odd:(3,2):d}) in degree $d$ is obtained using a
move $\widehat{\mu}_3$ (with $i=3$ and $j_3=1$) from the datum
with partitions
$$(3,2,\ldots,2),(d_{31}-1,d_{32},\ldots,d_{3(m_3-1)})$$
in degree $d-2$. This datum is constructible by the inductive
assumption, and the move is constructible because
$d_{31}-1\geqslant 2$, whence the conclusion.

\medskip

\noindent\textsc{Claim 4.} \emph{If
$(d_{2j})=(3,\ldots,3,2\ldots,2)$ then the datum can be realized
by a constructible diagram.} By induction on the odd number $k$ of
$3$'s in $(d_{2j})$. The base $k=1$ is given by Claim 3. Suppose
$k\geqslant 3$. Then it is easy to see that $(d_{3j})$ has at
least $k$ entries equal to $1$. Moreover by Claim 2 we can assume
$d_{31}\geqslant 3$. Then the datum~(\ref{S:S:odd:(3,2):d}) at
level $k$ in degree $d$ is obtained from the datum with partitions
$$(3,\ldots,3,2,\ldots,2),(d_{31}-1,d_{32},\ldots,d_{3(m_3-3)})$$
at level $k-2$ in degree $d-4$ by a move $\widehat{\mu}_1$ (with
$j_2=k-1$ and $j_3=1$) and a move $\widehat{\mu}_2$ (with $i=2$
and $j_2=m_2$) . The conclusion now follows from the inductive
assumption and from the remark that both moves are constructible,
because $\widehat{\mu}_1$ is applied to entries which are at least
$2$, and $\widehat{\mu}_2$ is applied to an entry $1$ created by
$\widehat{\mu}_1$.

To conclude the proof, consider a branch datum of
form~(\ref{S:S:odd:d}), and let us prove by induction on $d$ that
it is realized by a constructible diagram. The base step follows
from Fig.~\ref{degree5:fig}. To proceed, assume that no partition
is $(3,\ldots,3,2\ldots,2)$, which is not restrictive by Claim 4.
Claim 1 then implies that there are at least two $1$'s in
$(d_{ij})_{i=2,3}^{j=1,\ldots,m_i}$. If only $(d_{2j})$ contains
$1$'s then $d_{31}\geqslant 4$ by the assumption just made, so the
datum is obtained from a move $\widehat{\mu}_2$ (with $i=3$ and
$j_3=1$) from
$$(d_{21},\ldots,d_{2(m_2-2)}),(d_{31}-2,\ldots,d_{3m_3})$$
and the move is constructible because $d_{31}-2\geqslant 2$,
whence the conclusion.

Let us then assume that both $(d_{2j})$ and $(d_{3j})$ contain
some $1$. If one of these partitions, say $(d_{2j})$, is
$(2,\ldots,2,1)$, then we see that $d_{31}\geqslant 3$, and the
datum is obtained from
$$(2,\ldots,2,1),(d_{31}-1,\ldots,d_{3(m_3-1)})$$
by a move $\widehat{\mu}_3$ (with $i=3$ and $j_3=1$), which is
constructible because $d_{31}-1\geqslant 2$. If neither $(d_{2j})$
nor $(d_{3j})$ is $(2,\ldots,2,1)$, then the datum is obtained
from
$$(d_{21}-1,\ldots,d_{2(m_2-1)}),(d_{31}-1,\ldots,d_{3(m_3-1)})$$
by a move $\widehat{\mu}_1$ (with $j_2=j_3=1$), and the move is
constructible because either $d_{i1}-1\geqslant 2$ or there are
two $1$'s to choose from, one of which must be accessible by
Lemma~\ref{accessib:diag:lem}.\end{proof}

\paragraph{Coverings by the sphere: general case} Here we complete
the proof of Theorem~\ref{(d-2,2):sphere:thm}, using the results
and techniques of the previous paragraph. We begin with the following:

\begin{lemma}\label{refin:lem}
If $k\geqslant 1$ then the partitions $(y_j)_{j=1}^m$ of $2k$
with $m\geqslant k$ that do not refine the
partition $(k,k)$ are precisely the following:
\begin{itemize}
\item $(k+1,1,\ldots,1)$;
\item $(2,\ldots,2)$ when $k$ is odd.
\end{itemize}
\end{lemma}

\begin{proof}
We proceed by induction on $k$, the basis $k=1$ being obvious.
For the inductive step,
we assume as usual that $y_1\geqslant\ldots\geqslant y_m$. If $y_1>k$ then
the partition has the first special form listed and it does not refine $(k,k)$.
If $y_m\geqslant 2$ then the partition has the second special form listed and
it refines $(k,k)$ precisely when $k$ is even. We are left to show
that if $y_1\leqslant k$ and $y_m=1$ then $(y_j)_{j=1}^m$ refines $(k,k)$.
The conclusion is obvious if $y_1=1$, so we exclude this case.
Defining $x_1=y_1-1$ and $x_j=y_j$ for $j=2,\ldots,m-1$ we get a partition
$(x_j)_{j=1}^{m-1}$ of $2(k-1)$ to which the inductive assumption applies.
And it is actually very easy to check that $(y_j)$ refines $(k,k)$
both when $(x_j)$ refines $(k-1,k-1)$
and when $(x_j)$ has the second special form listed (it cannot have
the first special form).
\end{proof}

\dimo{(d-2,2):sphere:thm}
By Proposition~\ref{(d-2,2):sphere:odd:prop} it remains to
consider only the case of even degree. So we must prove realizability of
\begin{equation}\label{S:S:even}
\big(\matS,\matS,3,2k,(2k-2,2),(d_{2j}),(d_{3j})\big)
\end{equation}
with the two series of exceptions
\begin{eqnarray*}
& \big(\matS,\matS,3,2k,(2k-2,2),(2,\ldots,2),(2,\ldots,2)\big),\qquad k>2 & \\
& \big(\matS,\matS,3,2k,(2k-2,2),(2,\ldots,2),(k+1,1,\ldots,1)\big). &
\end{eqnarray*}
We first note that indeed
the data of these forms are non-realizable. For the former series
this is shown in~\cite[Corollary 6.4]{EKS}, while for the latter
this is a consequence of~\cite[Theorem 1.5]{partI}, because
$(d_{3j})$ does not refine $(k,k)$.

In the rest of the proof we exclude that~(\ref{S:S:even}) has one
of the exceptional forms and we show it is realizable. We first do
this assuming that all $d_{2j}$ are even. Since $d_{2j}\geqslant
2$, we have $m_2\leqslant k$ and hence $m_3\geqslant k$, so
Lemma~\ref{refin:lem} applies to $(d_{3j})$. And it allows us to
conclude that $(d_{3j})$ does refine $(k,k)$, because if it does
not the Riemann-Hurwitz condition implies that
$(d_{2j})=(2,\ldots,2)$, so we are in one of the cases excluded.

Let us then reorder $(d_{3j})$
so that $\sum_{j=1}^hd_{3j}=\sum_{j=h+1}^{m_3}d_{3j}=k$ for some $h$
and consider the branch datum
\begin{equation}\label{half:datum}
\big(\matS,\matS,4,k,(k-1,1),(d_{2j}/2),
(d_{3j})_{j=1}^h,(d_{3j})_{j=h+1}^{m_3}\big).
\end{equation}
We claim that it does not match any of the two data listed in
Proposition~\ref{(d-1,1):prop}. If (\ref{half:datum}) matches the first
datum of Proposition~\ref{(d-1,1):prop} then
$k=4$ and $m_2+m_3=6$, whereas we know that $m_2+m_3=2k$.
And (\ref{half:datum}) can only match the second datum of
Proposition~\ref{(d-1,1):prop}
if $(d_{3j})_{j=h+1}^{m_3}=(1,\ldots,1)$
up to permutation, whence $h=1$ by the same arguments used above,
and the datum is not matched anyway. We deduce that (\ref{half:datum})
is realizable, and the existence of a covering corresponding to~(\ref{S:S:even})
follows by composition with a covering realizing
$\big(\matS,\matS,3,2,(2),(2),(1,1)\big)$, which obviously exists.

\medskip

The argument carried out so far leaves us to consider the case where both
partitions $(d_{2j})$ and $(d_{3j})$ contain an odd entry. To deal with it
we follow the same line of reasoning as
in the proof of Proposition~\ref{(d-2,2):sphere:odd:prop}. To this
end we introduce a move $\mu$ on diagrams. The move depends on the
choice of a $G$-arc $e_i$ joining $\Gamma_{11}$ to some $\Gamma_{ij_i}$,
and its effect is described, depending on whether $i=2$ or $i=3$,
by either the top or the bottom part of Fig.~\ref{move1:fig}.
So $\mu_1$ is just the combination of two moves $\mu$ applied to
some $e_2$ and $e_3$.

Of course the effect of $\mu$ on the partitions $(d_{2j})$ and
$(d_{3j})$ is the move $\widehat{\mu}$ which consists in replacing
$d_{ij_i}$ by $d_{ij_i}+1$ and reordering, and appending $1$ to
the other partition. Not every move $\widehat{\mu}$ comes from
some $\mu$, but it does if $d_{ij_i}\geqslant 2$, since
in this case the tree $\Gamma_{ij_i}$ is accessible by
Lemma~\ref{accessib:diag:lem}.

\medskip

We now establish two claims that will readily lead to the conclusion.
Recall that both $(d_{2j})$ and $(d_{3j})$ have an odd entry

\medskip

\noindent\textsc{Claim 1.} \emph{At least one of $d_{2m_2}$
$d_{3m_3}$ is equal to 1.} Suppose the contrary. For $i=2,3$
not all $d_{ij}$ are even, so at least two of them are greater than 2
(recall that $(d_{ij})$ is a partition of $2k$), whence
$m_i\leqslant k-1$. But we know that $m_2+m_3=2k$, a contradiction.

\medskip

\noindent\textsc{Claim 2.} \emph{Up to switching
$(d_{2j})$ and $(d_{3j})$, we have $d_{21}\geqslant 3$ and
$d_{3m_3}=1$.}
By Claim~1 we can assume that $d_{3m_3}=1$. Suppose
that $d_{2j}\leqslant 2$ for all $j$. Since at least two of the
$d_{2j}$'s are odd, we must have $d_{2m_2}=d_{2(m_2-1)}=1$. In
particular $m_2\geqslant k+1$, hence $m_3\leqslant k-1$,
which implies that $d_{31}\geqslant 3$. Since $d_{2m_2}=1$
we get the desired situation by switching the partitions.

\medskip

To conclude the proof we remark that by Claim~2 the datum~(\ref{S:S:even})
can be obtained by the move $\widehat{\mu}$ from
$$\big(\matS,\matS,3,2k-1,(2k-3,2),(d_{21}-1,d_{22},\ldots,d_{2m_2}),(d_{31},\ldots,d_{3(m_3-1)})\big),$$
and this datum is realizable by some diagram thanks to
Proposition~\ref{(d-2,2):sphere:odd:prop}. In addition, the desired move is
constructible because $d_{21}-1\geqslant 2$. \finedimo

\paragraph{Coverings by the torus}
We now turn to the case $\Sigmatil=\matT$, where the result is even
stronger.

\begin{figure}
    \begin{center}
    \includegraphics[scale=0.5]{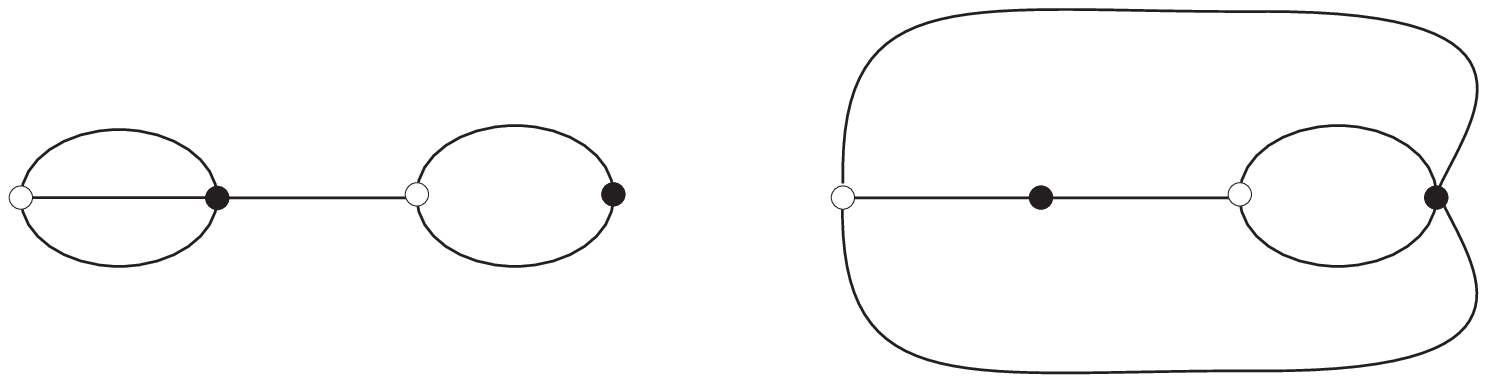}
    \mycap{Two abstract bipartite graphs.} \label{abstract4233:fig}
    \end{center}
    \end{figure}

\begin{figure}
    \begin{center}
    \includegraphics[scale=0.5]{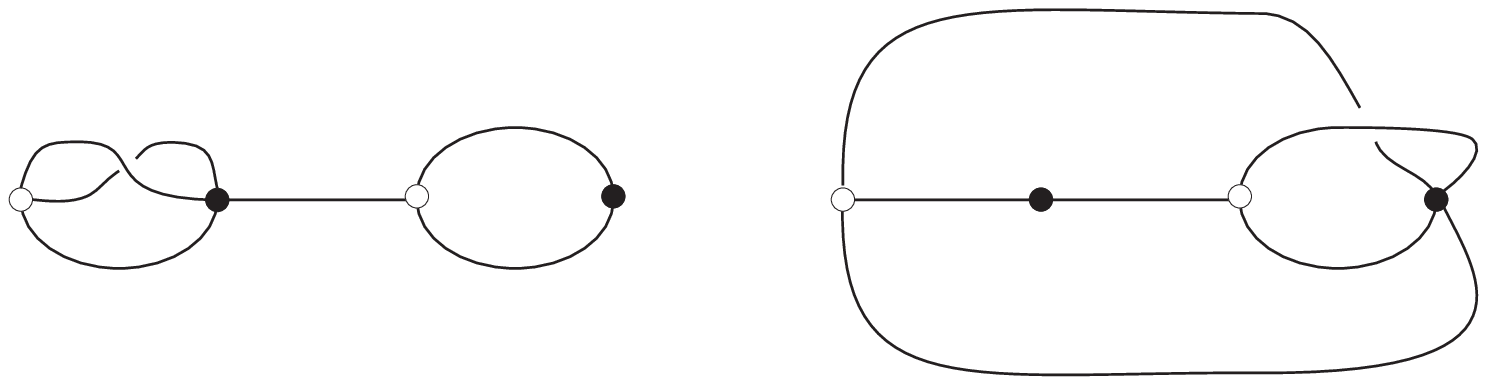}
    \mycap{Two dessins d'enfants on the torus.} \label{torus4233:fig}
    \end{center}
    \end{figure}

\dimo{(d-2,2):torus:thm} Note first that compatibility means
$m_2+m_3=d-2$. Let us first prove exceptionality of
$\big(\matT,\matS,3,6,(4,2),(3,3),(3,3)\big)$ using dessins
d'enfants, see Section~\ref{dessins:section}. There are only two
abstract bipartite graphs with two black vertices of valences 4
and 2 and two white vertices both of valence 3, see
Fig.~\ref{abstract4233:fig} (colors are used to represent the
bipartition).
It is now quite easy to analyze all the thickening of these graphs
giving a torus minus some discs. The number of discs is
automatically two, and one easily sees that one never gets two
discs of length 6, which would be necessary to realize the
partition $(3,3)$. For instance the thickenings given by taking
regular neighbourhoods of the planar immersions of
Fig.~\ref{torus4233:fig}
give the partitions $(5,1)$ and $(4,2)$ respectively. Note that
one could also use the second graph of Fig.~\ref{abstract4233:fig}
to get the partition $(5,1)$, but one could not the first graph of
Fig.~\ref{abstract4233:fig} to get the partition $(4,2)$.

Let us then prove realizability in general, excluding only the
case of $\big(\matT,\matS,3,6,(4,2),(3,3),(3,3)\big)$. To do this
we follow the same line of reasoning as in the proof of
Theorem~\ref{(d-2,2):sphere:thm}. Central to the proof is the move
$\mu$ introduced above. The following two facts concerning this
move are repeatedly used below:
\begin{itemize}
\item A diagram realizing a datum with $d_{ij}\geqslant 2$ for all
$i,j$ is accessible; \item $\mu$ transforms accessible diagrams
into accessible diagrams.
\end{itemize}
As in the previous proofs we proceed by establishing a series of claims.

\medskip

\noindent\textsc{Claim~1.} \emph{If $d$ is odd and
$d_{ij}\geqslant 2$ for all $i,j$ then $(d_{2j})$ is
$(3,2,\ldots,2)$ and $(d_{3j})$ is either $(5,2,\ldots,2)$ or
$(4,3,2,\ldots,2)$ or $(3,3,3,2,\ldots,2)$
up to permutation.} It follows from the assumptions that
$m_2,m_3\leqslant(d-1)/2$. Since their sum is $d-2$, the only
option is that one be $(d-1)/2$ and the other $(d-3)/2$. The
partitions of the statement of the claim are the only ones having
these lengths and not containing $1$'s. A similar argument proves
the next:

\medskip

\noindent\textsc{Claim~2.} \emph{If $d$ is even and
$d_{ij}\geqslant 2$ for all $i,j$ then one of the following
happens: \begin{itemize} \item Up to permutation,
$(d_{2j})=(2,\ldots,2)$ and $(d_{3j})\in\{(6,2,2,\ldots,2),\\
(5,3,2,\ldots,2),(4,4,2,\ldots,2), (4,3,3,2,\ldots,2),
(3,3,3,3,2,\ldots,2)\}$; \item $(d_{2j}),(d_{3j})\in
\{(4,2,\ldots,2),(3,3,2,\ldots,2)\}$.
\end{itemize}}

\medskip

\noindent\textsc{Claim~3.} \emph{If $d\geqslant 7$ and
$d_{ij}\geqslant 2$ for all $i,j$ then the datum is realizable}.
This is proved using the dessins d'enfant of
Section~\ref{dessins:section}.

    \begin{figure}
    \begin{center}
    \input{revised_new_dessins_on_torus.pstex_t} 
    \mycap{Dessins d'enfants on the torus for branch
    data not involving $1$.} \label{dessins_on_torus:fig}
    \end{center}
    \end{figure}
\begin{figure}
    \begin{center}
\includegraphics[scale=0.5]{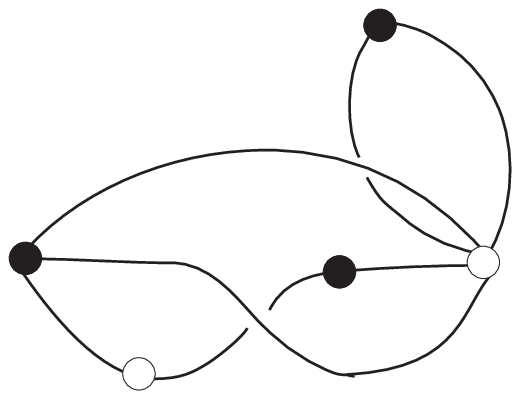}
    \mycap{A dessin d'enfant for the branch
    datum with partitions $(3,2,2),(4,3),(5,2)$.} \label{single_dessin_on_torus:fig}
    \end{center}
    \end{figure}

We present in Fig.~\ref{dessins_on_torus:fig} dessins d'enfants
proving realizability of all the desired data except that with
partitions $(3,2,2),(4,3),(5,2)$. For each dessin we have only
drawn the minimal number of vertices it should have: an infinite
series is obtained by adding pairs of vertices of opposite colors
on any of the edges not incident to the length-4 disc that one
easily sees in each picture. The branch datum with partitions
$(3,2,2),(4,3),(5,2)$ is realized in
Fig.~\ref{single_dessin_on_torus:fig}.

\medskip

\noindent\textsc{Claim~4.} \emph{If $d\leqslant 6$ then the datum
is realizable by an accessible diagram}. The first relevant $d$ is
$4$, where there is only the datum with partitions
$(2,2),(4),(4)$. This datum is realizable by
Theorem~\ref{full:cycle:thm}, and we conclude because
$d_{ij}\geqslant 2$ for all $i,j$. For $d=5$ we have $(3,2),
(3,2),(5)$ and $(3,2),(4,1),(5)$. Both are realizable by
Theorem~\ref{full:cycle:thm}, the former has $d_{ij}\geqslant 2$
for all $i,j$, and the latter is obtained from $(2,2),(4),(4)$ by
a move $\widehat{\mu}$, whence the conclusion. For $d=6$ we know
experimentally that only the excluded datum is exceptional. The
data with $d_{ij}\geqslant 2$ for all $i,j$ are realizable by
accessible diagrams, and those involving a $1$ are obtained via a
move $\widehat{\mu}$, and the conclusion follows from the case
$d=5$.

\medskip

\noindent\textsc{Claim~5.} \emph{If $d\geqslant 7$ then the datum
is realizable by an accessible diagram}. By induction. The base
step $d=7$ and the inductive step are essentially identical. For
both, either $d_{ij}\geqslant 2$ for all $i,j$, whence the
conclusion from Claim 3, or there is some $1$, so the datum is
obtained via a move $\widehat{\mu}$ from a datum in degree $d-1$.
For $d=7$ we use Claim 4 and the easy fact that the exceptional
datum can always be avoided, while for $d\geqslant 8$ the
conclusion is immediate from the inductive assumption.\finedimo

\paragraph{Conjunction of diagrams} In order to consider the case
of genus higher than 1, we describe a procedure that we
will use extensively.

Suppose we have two diagrams realizing branch data
$\big(\Sigmatil,\matS,3,d,(d_{ij})\big)$ and
$\big(\Sigmatil',\matS,3,d',(d_{ij}')\big)$ respectively. Let
$\alpha$ and $\alpha'$ be faces of opposite colors of these
diagrams. As was shown in the proof of Theorem~\ref{gTtoS:thm},
$\alpha$ is incident to three trees
$\Gamma_{1j_1},\Gamma_{2j_2},\Gamma_{3j_3}$ of distinct colors,
and, similarly, $\alpha'$ to some
$\Gamma'_{1j'_1},\Gamma'_{2j'_2},\Gamma'_{3j'_3}$. The following
steps define the \emph{conjunction of the diagrams along $\alpha$
and $\alpha'$}:
\begin{itemize}
\item Contract the
six trees $\Gamma_{ij_i}$ and $\Gamma'_{ij'_i}$
to distinct points and remove the interiors of $\alpha$ and
$\alpha'$. This gives two compact surfaces, both
bounded by a circle with three marked points bearing distinct colors 1,2,3.
\item Glue the two surfaces along their boundaries, matching the colors of
the marked points. This gives the closed surface
$\Sigmatil\#\Sigmatil'$ with an embedded checkerboard graph and
some trees.
\item Perform the white and black moves of
the proof of Theorem~\ref{gTtoS:thm} until a diagram with minimal
checkerboard graph is reached.
\end{itemize}
The resulting diagram of course satisfies the conditions of
Theorem~\ref{gTtoS:thm}, and one easily sees that it
realizes the datum
$\big(\Sigmatil\#\Sigmatil',\matS,3,d+d'-1, (d''_{ij})\big)$,
where for $i=1,2,3$ we get the partition
$(d''_{ij})$ by appending $(d'_{ij})$ to
$(d_{ij})$, removing the entries $d_{ij_i}$ and $d_{ij'_i}'$,
adding the entry $d_{ij_i}+d_{ij'_i}'-1$, and reordering. Note that
conjunction induces a certain move on branch data.
Constructibility of this move seems hard to investigate in general,
but we will prove it when necessary.

Let us note that the coloring of
faces and trees of one diagram can be changed arbitrarily,
so neither the condition that $\alpha$ and $\alpha'$ have opposite colors
nor the fact that the gluing be color-preserving is essential.
We also remark that the move $\mu$ used above can be viewed
as a conjunction with the diagram
realizing the datum $\big(\matS,\matS,3,2,(2),(2),(1,1)\big)$.

\paragraph{The case of higher genus} Here we consider the case
when the covering surface has genus at least 2. In this case there
are no obstructions to realizability.

\dimo{(d-2,2):high:genus:thm} We will actually establish the following
stronger statement: \emph{any compatible datum of the form
\begin{equation}\label{d-2:2:general:datum}
\big(g\matT,\matS,3,d,(d-2,2),(d_{2j})_{j=1}^{m_2},(d_{3j})_{j=1}^{m_3}\big),\qquad g\geqslant 2
\end{equation}
is realized by an accessible diagram}. Note that compatibility means that
$m_2+m_3=d-2g$, whence $d\geqslant 6$, and that the same statement
was already proved above for $g=1$, with a single exception for $d=6$.

The proof is by induction on $d$
simultaneously for all $g$. The base step is easy: for $d=6$,
there is a single compatible datum as in~(\ref{d-2:2:general:datum}), namely
$\big(2\matT,\matS,3,6,(4,2),(6),(6)\big)$, which is realizable by
Theorem~\ref{full:cycle:thm}. Then by Theorem~\ref{gTtoS:thm}
there is a diagram realizing the datum, and the diagram is automatically
accessible by Lemma~\ref{accessib:diag:lem}.

For the inductive step we need a preliminary fact. Recall that our partitions
are arranged non-increasingly.

\begin{lemma}\label{partitions:type:lem}
For a compatible datum as in~(\textrm{\ref{d-2:2:general:datum}}),
at least one of the following holds up to
switching $(d_{2j})$ and $(d_{3j})$:

\begin{enumerate}
    \item $d_{3m_3}=1$;
    \item All $d_{2j}$ and $d_{3j}$ are at least $2$, and $d_{3(m_3-1)}=d_{3m_3}=2$;
    \item $d_{2m_2}=d_{3m_3}=2$;
    \item There is an index $j$ such that $d_{3j}\geqslant 3$ and
    $d_{21}>d_{3j}$;
    \item All $d_{2j}$ and $d_{3j}$ are equal to one and the same integer $k$.
\end{enumerate}
\end{lemma}

\begin{proof}
Assume that neither of the first three cases holds, \emph{i.e.}
that $d_{2m_2}\geqslant 3$, $d_{3m_3}\geqslant 2$, and
$d_{3(m_3-1)}\geqslant 3$. If the
fourth case also does not hold, we deduce that $d_{21}\leqslant
d_{3(m_3-1)}$ and $d_{31}\leqslant d_{2m_2}$. Therefore
$$d_{21}\leqslant d_{3(m_3-1)}\leqslant d_{3(m_3-2)}\leqslant\ldots
\leqslant d_{31}\leqslant d_{2m_2}\leqslant
d_{2(m_2-1)}\leqslant\ldots d_{21},
$$
which implies that
$d_{21}=\ldots=d_{2m_2}=d_{31}=\ldots=d_{3(m_3-1)}$. Denoting this
number by $k$, we see that $d=km_2=k(m_3-1)+d_{3m_3}$. This
implies that $d_{3m_3}$ is divisible by $k$, and since
$0<d_{3m_3}\leqslant d_{3(m_3-1)}=k$, it is also equal to $k$.
Thus, we have the fifth case.
\end{proof}

We now proceed with the inductive step, supposing $d>6$
and the realizability of~(\ref{d-2:2:general:datum}) to be known
for degrees smaller than $d$. We
treat the different cases described in
Lemma~\ref{partitions:type:lem} separately, employing both the technique of
dessins d'enfants and the criterion of Theorem~\ref{gTtoS:thm}.
We also notice
that by Lemma~\ref{accessib:diag:lem} a realizable datum that does
not involve 1's is always realized by an accessible diagram, so in
Cases~2-5 (where we always assume that Case~1 is not applicable)
it will be sufficient to show that our datum is realizable.
For this reason, recalling Theorem~\ref{full:cycle:thm},
in Cases~2-5 we will also assume $m_2,m_3\geqslant 2$.

\medskip

\noindent\textsc{Case~1:} $d_{3m_3}=1$.
In this case~(\ref{d-2:2:general:datum})
is obtained by the move $\widehat{\mu}$ (described in
the proof of Theorem~\ref{(d-2,2):sphere:thm}) from the datum
$$\big(g\matT,\matS,3,d-1,(d-2,2),
(d_{21}-1,d_{22},\ldots,d_{2m_2}),(d_{31},\ldots,d_{3(m_3-1)})\big).$$
This move is constructible, since the latter datum can be realized
by an accessible diagram by the inductive assumption. The arising
diagram is also accessible, since the move $\mu$ preserves accessibility.

\medskip

\noindent\textsc{Case~2:}
\emph{$d_{2m_2}\geqslant 2$ and $d_{3m_3}=d_{3(m_3-1)}=2$}.
We first note that we can assume $d_{21}>2$, because
if $d_{21}=2$ then either $d_{31}>2$, and we can switch the
partitions, or $d_{ij}=2$ for $i=2,3$ and all $j$, which is only
possible for $g=0$. Now we consider the compatible datum
$$\big(g\matT,\matS,3,d-4,(d-4),(d_{21}+d_{22}-4,d_{23},\ldots,d_{2m_2}),
(d_{31},\ldots,d_{3(m_3-2)})\big).$$ By
Theorem~\ref{full:cycle:thm} it is realizable, so there exists
in $g\matT$ a dessin d'enfant $D$ realizing it and such that
$g\matT\setminus D$ is a disc of length $2(d-4)$. We take the
vertex of $D$ corresponding to the entry $d_{21}+d_{22}-4$ and we
perform at it the move shown in
Fig.~\ref{dessin:(2,2):fig}.
\begin{figure}
    \begin{center}
    \includegraphics[scale=0.5]{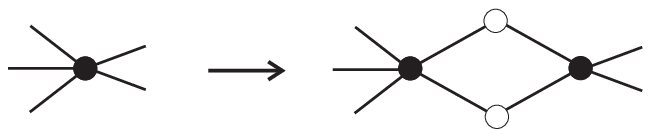}
    \mycap{A move on dessins d'enfants.}\label{dessin:(2,2):fig}
    \end{center}
    \end{figure}
Clearly, we can do so in such a way that the new black
vertices have valences $d_{21}$ and $d_{22}$. Hence we get a
dessin d'enfant realizing~(\ref{d-2:2:general:datum}).

\medskip

\noindent\textsc{Case~3:}
\emph{$d_{2m_2}=d_{3m_3}=2$.} Since $m_2,m_3\geqslant 2$ we
easily see that $d\geqslant 8$, so the inductive assumption implies realizability
of the compatible datum
$$\big(g\matT,\matS,3,d-2,(d-4,2),(d_{21},\ldots,d_{2(m_2-1)}),
(d_{31},\ldots,d_{3(m_3-1)})\big).$$ Let $D$ be a dessin d'enfant
in $g\matT$ realizing it and such that $g\matT\setminus D$
consists of two discs of lengths $2(d-4)$ and $4$. Since
$2(d-4)>4$ there exists an edge $e$ of $D$ to which the first disc
is doubly incident. We then place on $e$ two vertices of valence 2
and give them colors so to get a bipartite graph. The result is of
course a dessin d'enfant realizing~(\ref{d-2:2:general:datum}).

\medskip

\noindent\textsc{Case~4:} \emph{$d_{2m_3}\geqslant 3$,
$d_{3(m_3-1)}\geqslant 3$, $d_{3m_3}\geqslant 2$ and
$d_{21}>d_{3j}\geqslant 3$ for some $j$.} This implies that either
$d_{21}>d_{3m_3}\geqslant 3$ or $d_{21}>d_{3(m_3-1)}$ and
$d_{3m_3}=2$. Let $k$ be the smallest entry among $d_{2j}$,
$d_{3j}$ that is not equal to $2$. Clearly, $k$ can only be
$d_{3(m_3-1)}$, $d_{3m_3}$, or $d_{2m_2}$. The proof is basically
the same in all three cases. We then assume the first one occurs
(which implies $d_{3m_3}=2$), but later we will mention the
variations for the other two cases. The current assumption
$m_2,m_3\geqslant 2$ easily implies that $k\leqslant d-3$.

\medskip

\noindent\textsc{Case~4.1:}
\emph{$k$ is odd.}
We consider the compatible data
\begin{equation}\label{case:4:1:auxbis}
\begin{array}{rr}
\big(\gamma\matT,\matS,3,d-k,(d-2-k,2),
(d_{21}-k,d_{22},\ldots,d_{2m_2}),\ \ \\
(d_{31},\ldots,d_{3(m_3-2)},d_{3m_3})\big)
\end{array}
\end{equation}
with $\gamma=g-(k-1)/2=(d-k-(m_2+m_3)+1)/2$, and
\begin{equation}\label{case:4:1:aux}
\big((k-1)/2\cdot\matT,\matS,3,k+1,(k+1),(k+1),(k,1)\big).
\end{equation}
We first show that
if~(\ref{case:4:1:auxbis}) is non-realizable then~(\ref{d-2:2:general:datum})
is realizable. By the inductive assumption, Proposition~\ref{(d-2,2):sphere:odd:prop} and
Theorem~\ref{(d-2,2):torus:thm} we know~(\ref{case:4:1:auxbis}) can only be
non-realizable if $\gamma=0$ or~(\ref{case:4:1:auxbis}) is
$\big(\matT,\matS,3,6,(4,2),(3,3),(3,3)\big)$. In the second case~(\ref{d-2:2:general:datum})
has the form
$$\big(g\matT,\matS,3,k+6,(k+4,2),(k+3,3),(3,k,3)\big)$$
but by the definition of $k$ we must have $k=3$, so~(\ref{d-2:2:general:datum}) is
$$\big(2\matT,\matS,3,9,(7,2),(6,3),(3,3,3)\big),$$
which is realized in Fig.~\ref{dessin:(7,2):fig}.
\begin{figure}
    \begin{center}
    \includegraphics[scale=0.5]{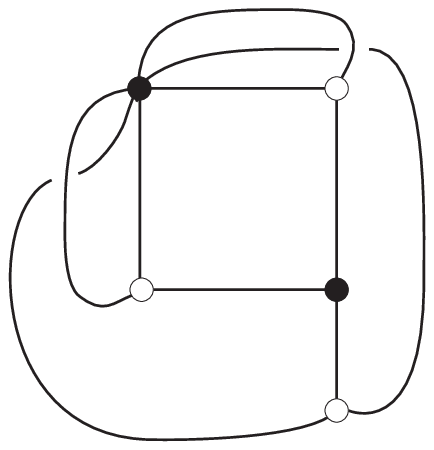}
    \mycap{A dessin d'enfant realizing
    $\big(2\matT,\matS,3,9,(7,2),(6,3),(3,3,3)\big)$.} \label{dessin:(7,2):fig}
    \end{center}
    \end{figure}
Suppose $\gamma=0$. Since $m_2\leqslant d/k$ and $m_3\leqslant (d-2)/k+1$, we have
$$0=2\gamma\geqslant d-k-2(d-1)/k.$$
Knowing that $3\leqslant k\leqslant d-3$, a tiny bit of algebra
shows that this can only happen for $d=7$ (and $k=3,4$), but we
are assuming $m_2,m_3\geqslant 2$ and $g\geqslant 2$, so the
compatibility condition $m_2+m_3=d-2g$ is impossible for $d=7$.

We can then assume~(\ref{case:4:1:auxbis}) is realizable, so it is
realized by a diagram which is automatically accessible. Within
any such diagram we can select a face $\alpha$ incident to $\Gamma_{11}''$ and
$\Gamma_{21}''$, whence of type $(d-2-k,d_{21}-k,d_{3j})$
for some $j$. By Theorem~\ref{full:cycle:thm} the datum~(\ref{case:4:1:aux})
is also realizable by a diagram, and of course we can find in the diagram
a face $\alpha'$ of type $(k+1,k+1,1)$. Performing the
conjunction along $\alpha$ and $\alpha'$
of the diagrams realizing~(\ref{case:4:1:auxbis}) and~(\ref{case:4:1:aux})
we get a diagram realizing~(\ref{d-2:2:general:datum}).

\medskip

\noindent\textsc{Case~4.2:}
\emph{$k$ is even.}
We consider the compatible data
\begin{equation}\label{case:4:2:auxbis}
\big(\gamma\matT,\matS,3,d-k,(d-k),
(d_{21}-k,d_{22},\ldots,d_{2m_2}),
(d_{31},\ldots,d_{3(m_3-2)},d_{3m_3})\big)
\end{equation}
with $\gamma=g-(k-2)/2=(d-k-(m_2+m_3)+2)/2$, and
\begin{equation}\label{case:4:2:aux}
\big((k-2)/2\cdot\matT,\matS,3,k+1,(k-1,2),(k+1),(k,1)\big).
\end{equation}
Both these data are realizable by Theorem~\ref{full:cycle:thm}.
Any diagram realizing~(\ref{case:4:2:auxbis}) has a face $\alpha$
of type $(d-k,d_{21}-k,d_{3j})$ for some $j$. Moreover, since $k>2$,
the inductive assumption implies that~(\ref{case:4:2:aux})
can be realized by an accessible diagram, so we can choose in it
a face $\alpha'$ incident to trees of degrees $k-1$ and $1$.
Therefore $\alpha'$ has type $(k-1,k+1,1)$. Just as above, we
do the conjunction along $\alpha$ and $\alpha'$
of the diagrams realizing~(\ref{case:4:2:auxbis}) and~(\ref{case:4:2:aux}),
getting a diagram realizing~(\ref{d-2:2:general:datum}).

\medskip

Recall now that in Cases~4.1-2
we have assumed that $k$ is $d_{3(m_3-1)}$, whereas it can also
be $d_{3m_3}$ or $d_{2m_2}$, but the arguments given readily extend to these
cases. If $k$ is $d_{3m_3}$ we replace the second and third partitions
of $d$ in both~(\ref{case:4:1:auxbis}) and~(\ref{case:4:2:auxbis}) by
$$(d_{21}-k,d_{22},\ldots,d_{2m_2}),\quad
(d_{31},\ldots,d_{3(m_3-1)}).$$
If $k$ is $d_{2m_2}$ we replace them by
$$(d_{21},\ldots,d_{2(m_2-1)}),\quad
(d_{31}-k,d_{32},\ldots,d_{3m_3}).$$

\medskip

\begin{figure}
    \begin{center}
    \includegraphics[scale=0.5]{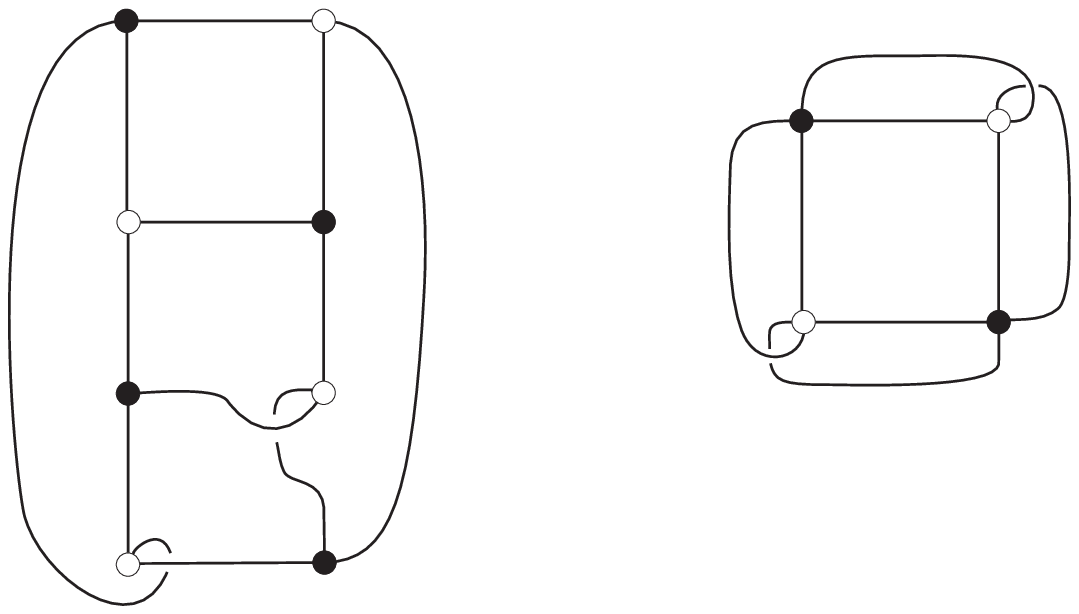}
    \mycap{Dessins d'enfants realizing the data
    $\big(2\matT,\matS,3,12,(10,2),(3,3,3,3),(3,3,3,3)\big)$
    and $\big(2\matT,\matS,3,8,(6,2),(4,4),(4,4)\big)$.} \label{dessins:3344:fig}
    \end{center}
    \end{figure}

\noindent\textsc{Case~5:}
\emph{all $d_{2j}$ and $d_{3j}$ are equal to some $k$.}
We must prove realizability of a compatible datum of the form
\begin{equation}\label{general:kh}
\big(g\matT,\matS,3,kh,(kh-2,2),(k,\ldots,k),(k,\ldots,k)\big),\qquad
g\geqslant 2,\ k\geqslant 3,\ h\geqslant 1.
\end{equation}
Compatibility means that $kh=2(g+h)$, so either $k$ or $h$ is even,
and one easily sees that one of the following happens:
\begin{itemize}
\item[(a)] $k=3$ and $h\geqslant 4$ is even;
\item[(b)] $k=4$ and $h\geqslant 2$;
\item[(c)] $k\geqslant 6$ is even and $h\geqslant 1$;
\item[(d)] $k\geqslant 5$ is odd and $h\geqslant 2$ is even.
\end{itemize}
In cases (a,b,c) we will prove realizability of~(\ref{general:kh}) fixing
$k$ and proceeding by induction on $h$, with an induction step of length 2 in
case (a). In case (d) we will give a direct (but very similar) argument.
The base step of the induction in cases (a) and (b) is established in
Fig.~\ref{dessins:3344:fig},
and it is a consequence of Theorem~\ref{full:cycle:thm} in case (c). The core
of our arguments is the following:

\medskip

\noindent\textsc{Claim.} \emph{Suppose that~(\ref{general:kh})
is realizable. Then
\begin{equation}\label{claim:aux}
\big(g\matT,\matS,3,kh+2,(kh,2),(k,\ldots,k,2),\\
(k,\ldots,k,2)\big)
\end{equation}
can be realized by a diagram having a face of type $(kh,2,2)$.}
Indeed, consider a dessin d'enfant $D\subset\Sigmatil$
realizing~(\ref{general:kh}) and such that all vertices of $D$ have valence
$k$. Then the complement of $D$ consists of two
discs $B'$ and $B''$, of lengths $2(kh-2)$ and $4$ respectively.
Since $kh\geqslant 6$ we have $2(kh-2)>4$, so there is an edge $e$ of $D$ to which $B'$
is doubly incident. Then we place two
vertices $u$ and $v$ of valence 2 on $e$ and, as in the proof of Case 3
above, we get a dessin d'enfant $\widetilde{D}$ realizing~(\ref{claim:aux}).

It follows from the proof of Theorem~\ref{gTtoS:thm} that a
diagram realizing~(\ref{claim:aux}) can now be
constructed as follows. Choose a point $x'$ in the interior of
$B'$ and take a cone (in $B'$) with center $x'$ over the pull-back
of the vertex set of $\widetilde{D}$ in the boundary of $B'$. Do
the same for some interior point $x''$ of $B''$. We get a
triangulation of $\Sigmatil$ that we can color in a
checkerboard fashion,
and we notice that there is a
triangle incident to both $u$ and $v$ (two of them, actually, a
white one and a black one). Assign color 1 to $x'$ and
$x''$, and colors 2 and 3 to the members of the two different
partitions of the vertex set of the bipartite graph
$\widetilde{D}$. It is easy to see now that applying to this
triangulation as long as possible the black and white moves
described in the proof of
Theorem~\ref{gTtoS:thm} yields a diagram with
a white and a black disc of type $(kh,2,2)$. The Claim is proved.

\medskip

\noindent\emph{Inductive step in case} (a). The inductive assumption is that
some datum
\begin{equation*}
\big(h/2\cdot\matT,\matS,3,3h,(3h-2,2),(3,\ldots,3),(3,\ldots,3)\big)
\end{equation*}
is realizable. By the Claim we can then realize
\begin{equation}\label{kh:k=3:aux}
\big(h/2\cdot\matT,\matS,3,3h+2,(3h,2),(3,\ldots,3,2),(3,\ldots,3,2)\big)
\end{equation}
via a diagram with a face $\alpha$ of type $(3h,2,2)$. Now, starting from the obvious
dessin d'enfant realizing
\begin{equation}\label{kh:k=3:auxbis}
\big(\matT,\matS,3,5,(5),(3,2),(3,2)\big)
\end{equation}
and applying the triangulation trick used to prove the Claim,
we see that this datum is realized by a diagram with a face $\alpha'$
of type $(5,2,2)$. Taking the conjunction of the diagrams
realizing~(\ref{kh:k=3:aux}) and~(\ref{kh:k=3:auxbis})
along $\alpha$ and $\alpha'$ we get a realization of the datum
\begin{equation*}
\big((h+2)/2\cdot\matT,\matS,3,3(h+2),(3h+4,2),(3,\ldots,3),(3,\ldots,3)\big)
\end{equation*}
\emph{i.e.}~the desired conclusion.

\medskip

\noindent\emph{Inductive step in cases} (b,c).
The inductive assumption is that
some datum
\begin{equation*}
\big(g\matT,\matS,3,kh,(kh-2,2),(k,\ldots,k),(k,\ldots,k)\big)
\end{equation*}
with $g=(k-2)h/2$ is realizable. By the Claim we can then realize
\begin{equation}\label{kh:even_k:aux}
\big(g\matT,\matS,3,kh+2,(kh,2),(3,\ldots,3,2),(3,\ldots,3,2)\big)
\end{equation}
via a diagram with a face $\alpha$ of type $(kh,2,2)$. Now note
that the datum
\begin{equation}\label{kh:even_k:auxbis}
\big((k-2)/2\cdot\matT,\matS,3,k-1,(k-1),(k-1),(k-1)\big)
\end{equation}
is compatible, whence realizable by Theorem~\ref{full:cycle:thm}.
Taking the conjunction of the diagram
realizing~(\ref{kh:even_k:aux}) and any diagram
realizing~(\ref{kh:even_k:auxbis}) along $\alpha$ and any face of
the latter we get a realization of the datum
\begin{equation*}
\big(\gamma\matT,\matS,3,k(h+1),(k(h+1)-2,2),(k,\ldots,k),(k,\ldots,k)\big)
\end{equation*}
with $\gamma=(k-2)(h+1)/2$. This is what we had to prove.

\medskip

\noindent\emph{Proof of realizability in case} (d). We begin by
noting that if $h\geqslant 2$ is even and $k\geqslant 5$ is odd
then the datum
\begin{equation*}
\big(\gamma\matT,\matS,3,k(h-1),(k(h-1)),(k,\ldots,k),(k,\ldots,k)\big),
\end{equation*}
with $\gamma=((k-2)(h-1)+1)/2$, is compatible, whence realizable
by Theorem~\ref{full:cycle:thm}. The very same argument used to
prove the Claim then shows that the datum
\begin{equation}\label{kh:odd_k:aux}
\big(\gamma\matT,\matS,3,k(h-1)+2,(k(h-1)+2),(k,\ldots,k,2),(k,\ldots,k,2)\big)
\end{equation}
is realized by a diagram with a face $\alpha$ of type
$(k(h-1)+2,2,2)$. Using the same techniques one can also easily
show that the datum
\begin{equation}\label{kh:odd_k:auxbis}
\big((k-3)/2\cdot\matT,\matS,3,k-1,(k-3,2),(k-1),(k-1)\big)
\end{equation}
is realized by a diagram with a face $\alpha'$ of type $(k-3,k-1,k-1)$.
And again the conclusion follows because the conjunction of
the diagrams realizing~(\ref{kh:odd_k:aux}) and~(\ref{kh:odd_k:auxbis})
along $\alpha$ and $\alpha'$ gives a
realization of~(\ref{general:kh}) under the assumptions of (d).

\medskip

Having treated all five cases of Lemma~\ref{partitions:type:lem} we have eventually shown
Theorem~\ref{(d-2,2):high:genus:thm}.\finedimo

\vspace{.35cm}

\noindent 
Chelyabinsk State University\\
ul. Br. Kashirinykh, 129\\
454021 Chelyabinsk, Russia\\
pervova@csu.ru

\vspace{.35cm}

\noindent 
Dipartimento di Matematica Applicata\\
Universit\`a di Pisa\\
Via Bonanno Pisano 25B\\
56126 Pisa, Italy\\
petronio@dm.unipi.it

\end{document}